\documentclass[11pt]{article}
\usepackage[latin1]{inputenc}
\usepackage[english]{babel}
\usepackage{amsmath}
\usepackage{paralist}
\usepackage{amssymb,amsfonts, color}
\usepackage[margin=2.7cm]{geometry}
\usepackage{xcolor}
\usepackage{graphicx, color, enumerate}
\usepackage[latin1]{inputenc}
\usepackage[active]{srcltx}
\usepackage{tikz}
\usepackage{pgf}
\usepackage{tkz-euclide}
\usepackage{etex}
\usepackage{verbatim}
\usepackage{tikz-3dplot}
\usepackage{pgfkeys}
\usepackage{fp}

\newtheorem{theorem}{Theorem}[section]
\newtheorem{proposition}[theorem]{Proposition}
\newtheorem{corollary}[theorem]{Corollary}
\newtheorem{lemma}[theorem]{Lemma}
\newtheorem{remark}[theorem]{Remark}
\newtheorem{definition}[theorem]{Definition}
\newtheorem{example}[theorem]{Example}

\newcommand{\bcl}{\begin{center}}
\newcommand{\ecl}{\end{center}}
\newcommand{\brl}{\begin{right}}
\newcommand{\erl}{\end{right}}
\newcommand{\ben}{\begin{enumerate}}
\newcommand{\een}{\end{enumerate}}
\newcommand{\overliner}{\begin{array}}
\newcommand{\earr}{\end{array}}
\newcommand{\btab}{\begin{tabular}}
\newcommand{\etab}{\end{tabular}}
\newcommand{\bdoc}{\begin{document}}
\newcommand{\edoc}{\end{document}}
\newcommand{\beqy}{\begin{eqnarray}}
\newcommand{\eeqy}{\end{eqnarray}}

\newcommand{\beqi}{\begin{eqnarray*}}
\newcommand{\eeqi}{\end{eqnarray*}}
\newcommand{\bitem}{\begin{itemize}}

\newcommand{\eitem}{\end{itemize}}
\newcommand{\nln}{\newline}
\newcommand{\newt}{\newtheorem}

\newcommand{\pa}{\partial}
\newcommand{\re}{{I\!\!R}}
\newcommand{\ren}{\re^N}
\newcommand{\xr}{x\in\re }
\newcommand{\x}{\times}
\newcommand{\dyle}{\displaystyle}
\newcommand{\ene}{{I\!\!N}}
\newcommand{\irn}{\int\limits_{\re^N}}
\newcommand{\io}{\int\limits_{\O}}
\newcommand{\meas}{{\rm meas\,}}

\newcommand{\sign}{{\rm sign}}
\newcommand{\map}{\longrightarrow }
\newcommand{\imp}{\Longrightarrow }
\renewcommand{\div}{\nabla\cdot }
\newcommand{\sen}{{\rm sen\,}}
\newcommand{\tg}{{\rm tg\,}}
\newcommand{\arcsen}{{\rm arcsen\,}}
\newcommand{\arctg}{{\rm arctg\,}}
\newcommand{\supp}{{\textsl supp\ }}
\newcommand{\ity}{\int_{-\iy}^{+\iy}}
\newcommand{\limit}{\lim\limits}
\newcommand{\limi}{\limit_{n\to\infty}}
\newcommand{\sumi}{\sum\limits_{n=1}^{\infty}}
\newcommand{\ulu}{\underline u}
\newcommand{\ulw}{\underline w}
\newcommand{\ulz}{\underline z}
\newcommand{\ulv}{\underline v}
\newcommand{\uls}{\underline s}
\newcommand{\olu}{\overline u}
\newcommand{\olv}{\overline v}
\newcommand{\ols}{\overline s}
\newcommand{\ob}{\overline\b}
\newcommand{\ovar}{\overline\var}
\newcommand{\wv}{\widetilde v}
\newcommand{\wu}{\widetilde u}
\newcommand{\ws}{\widetilde s}
\renewcommand{\a }{\alpha }
\renewcommand{\b }{\beta }
\newcommand{\g }{\gamma}
\newcommand{\G }{\Gamma }
\renewcommand{\d }{\delta }

\newcommand{\D }{\Delta }
\newcommand{\e }{\varepsilon }
\newcommand{\z }{\zeta }
\renewcommand{\l }{\lambda }
\renewcommand{\L }{\Lambda }
\newcommand{\m }{\mu }
\newcommand{\n }{\nu }
\newcommand{\s }{\sigma }
\renewcommand{\r }{\rho }
\newcommand{\Sig }{\Sigma }
\renewcommand{\t }{\tau }
\newcommand{\var }{\varphi }
\renewcommand{\o }{\omega }
\renewcommand{\O }{\Omega }
\newcommand{\bR}{{\bf R}}
\newcommand{\bC}{{\bf C}}
\newcommand{\bZ}{{\bf Z}}
\newcommand{\bN}{{\bf N}}
\newcommand{\bQ}{{\bf Q}}
\newcommand{\bK}{{\bf K}}
\newcommand{\bI}{{\bf I}}
\newcommand{\bv}{{\bf v}}
\newcommand{\bV}{{\bf V}}
\DeclareMathOperator{\suppo}{supp} \DeclareMathOperator{\di}{div}

\def\qed{\unskip\kern 6pt \penalty 500
\raise -2pt\hbox{\vrule \vbox to10pt{\hrule width 4pt \vfill\hrule}\vrule}\par}

\newenvironment{Proof}{\removelastskip\vskip12pt
plus 1pt \noindent\em\rm}{\hfill {\qed \hskip .2cm}}

\title{Relaxed uniqueness conditions for the\\ parabolic Schr\"odinger equation \\ on Riemannian manifolds}

\author{Fabio
Punzo\thanks{Dipartimento di Matematica, Politecnico di Milano, Via Bonardi 9, I-20133 Milano, Italy
(fabio.punzo@polimi.it). }
}

\date{}

\begin{document}
\maketitle

\abstract{We study uniqueness for solutions to the Cauchy problem associated with
the parabolic Schr\"odinger equation on complete noncompact Riemannian manifolds,
under suitable integral conditions on the solution. We show that, under suitable assumptions on the potential
$V$, the required integrability condition can be significantly relaxed compared to the case without potential. This improvement is achieved by exploiting the decay of positive solutions to the associated stationary Schr\"odinger equation. To the best of our knowledge, identifying how the behavior of the potential influences the uniqueness integral condition, through the decay properties of solutions to the corresponding stationary equation, constitutes a novel contribution to the theory.}

\medskip
\smallskip

\bigskip
\smallskip

\section{Introduction}\setcounter{equation}{0}
We investigate the uniqueness of solutions, under certain integral conditions, to the initial value problem
\begin{equation}\label{e1}
\begin{cases}
u_t - \Delta u + V u = 0 & \text{in } M\times (0, T] \\
u = 0 & \text{in } M\times \{0\},
\end{cases}
\end{equation}
where $(M,g)$ is a complete noncompact Riemannian manifold, $T>0$, $\Delta$ denotes the Laplace-Beltrami operator, and $V\in C^1(M)$.

The uniqueness of solutions to linear second-order parabolic equations posed on $\mathbb R^n$ has been extensively studied (see, e.g., \cite{ArBe, EKP, KP, PoPuTe1, PoPuTe2, Sm, Son}). This question has also been addressed in the context of Riemannian manifolds. In particular, the uniqueness of bounded solutions to the problem
\begin{equation}\label{e1h}
\begin{cases}
u_t - \Delta u  = 0 & \text{in } M\times (0, T]\\
u = 0 & \text{in } M\times \{0\}
\end{cases}
\end{equation}
has been the subject of extensive investigation. A key result in this direction is the equivalence between uniqueness and the \emph{stochastic completeness} of the manifold \( M \); see \cite{Grig1}. Stochastic completeness is characterized by the condition
\[
\int_M p(x, y, t) \, d\nu(y) = 1, \quad \text{for all } x \in M, \; t > 0,
\]
where \( p(x, y, t) \) is the heat kernel and \( d\nu \) is the Riemannian volume element.

According to \cite{Yau}, any geodesically complete Riemannian manifold with Ricci curvature bounded from below is stochastically complete. This result was extended in \cite{Ich, hsu1989heat} to allow for negative Ricci curvature under appropriate lower bounds. Under similar geometric assumptions, the uniqueness of \( L^1(M) \)-solutions was established in \cite{li1984uniqueness}.

To explore further results, fix \( x_0 \in M \) and define the geodesic ball
\[ B_R(x_0) := \{ x \in M : r(x) < R\}, \]
where \( r(x) := \operatorname{dist}(x, x_0) \). Let \( \mathcal V(x_0, R) \) denote the volume of \( B_R(x_0) \). In \cite{gaffney1959conservation}, it is shown that if
\[
\log \mathcal V(x_0, R) = o(R) \quad \text{as } R \to \infty,
\]
then \( M \) is stochastically complete. Another sufficient condition due to \cite{Grig0} states that if
\[
\int_1^\infty \frac{r}{\log \mathcal V(x_0, r)} \, dr = \infty,
\]
and \( M \) is geodesically complete, then it is stochastically complete. We refer the reader to \cite{Grig1} for a comprehensive overview (see also \cite{Mu1}, \cite{Mu2}).

It is proved in \cite[Theorem 9.2]{Grig1} for $p=2$ and in \cite[Theorem 2.2]{Pu15} for $1 \leq p \leq 2$ that if $u$ is a solution to \eqref{e1h}, and there exist \( x_0 \in M \), \( R_0 > 0 \), and an increasing continuous function \( h: (0,\infty) \to (0,\infty) \) such that
\[
\int_{R_0}^\infty \frac{r}{\varphi(h)} \, dr = \infty,
\]
and
\begin{equation}\label{ei5}
\int_0^T \int_{B_R(x_0)} |u(x,t)|^p \, d\nu(x) dt \leq \exp\{h(R)\} \quad \text{for all } R > R_0,
\end{equation}
then \( u \equiv 0 \) on \( M \times (0,T] \). In particular, this implies that if
\begin{equation}\label{ei6}
\int_0^T \int_M |u(x,t)|^p \, \exp\{-h(r(x))\} \, d\nu(x) dt < \infty,
\end{equation}
then \( u \equiv 0 \) on \( M \times (0,T] \).

\smallskip

Similar uniqueness results hold for the heat equation on a weighted Riemannian manifold $(M,g,\mu)$, where the measure $\mu$ is defined through a positive, sufficiently smooth function $\varphi : M \to \mathbb{R}$ by
$d\mu = \varphi^2 \, d\nu.$ In this context, the weighted Laplacian $\Delta_\mu$ is used (see Section \ref{mb} for details). In \cite[Theorem 11.9]{Grig3} for $p=2$ and \cite{MeRo} for $p\neq 2$ (see also \cite[Theorem 2.2]{Pu15}), it is shown that the solution $u \equiv 0$ is unique for the problem
\begin{equation}\label{e35}
\begin{cases}
  w_t - \Delta_{\mu} w = 0 & \text{in } M\times (0, T] \\
  w = 0 & \text{in } M\times \{0\},
\end{cases}
\end{equation}
provided that
\begin{equation}\label{ei7}
\int_0^T \int_M |u(x,t)|^p \, \exp\{-h(r(x))\} \, d\mu(x) dt < \infty.
\end{equation}
\medskip

We now discuss results concerning problem \eqref{e1} with a potential term \( V \). These are presented in \cite{IM}, where more general operators are also considered. Let \( h \) be as above and set \( V_- := \max\{-V, 0\} \). It is shown in \cite[Theorem 2.1]{IM} that if \( u \) solves \eqref{e1} and satisfies \eqref{ei5} with \( p=2 \), then \( u \equiv 0 \), provided that, for some \( C > 0 \),
\[
\int_{B_R(x_0)} V_-(x) \psi^2(x) \, d\nu \leq C\left[\frac{h(R)}{R}\right]^2 \int_{B_R(x_0)} \psi^2(x) \, d\nu
\]
for any \( R > 1 \) and \( \psi \in C_c^\infty(B_R(x_0)) \). Thus, the integral condition \eqref{ei5} is still requested despite the presence of the potential \( V \).

\medskip

Let us point out that in \cite[Section 10]{Grig2} and \cite{GSC}, estimates for  the heat kernel \( p^V \) associated with the Schr\"odinger operator
\( -\Delta + V \) are derived.  It is observed that, in some cases, \( p^V \) behaves similarly to \( p \),
the heat kernel of the Laplacian \( -\Delta \); however, under suitable  conditions on the potential \( V \), the behavior of \( p^V \) is  significantly influenced by \( V \).

Since, broadly speaking, the uniqueness class is generally related to the behavior of the heat kernel, the aforementioned results suggest that the uniqueness class for equation \eqref{e1} may depend on the potential $V$. In fact, the aim of this paper is to show that, under suitable assumptions on the potential \( V \), the integral condition on the solution, ensuring uniqueness, can be relaxed. This refinement stems from an analysis of the decay properties of positive solutions to the corresponding stationary Schr\"odinger equation. As far as we are aware, the idea that the uniqueness integrability condition can be relaxed in relation to the asymptotic behavior of positive solutions to the stationary equation, appears to be a novel perspective in the existing literature.

Specifically, instead of requiring \eqref{ei6}, we impose a much weaker requirement, related to a positive solution of the stationary Schr\"odinger equation
\begin{equation}\label{ei2}
\Delta \varphi - V(x)\varphi = 0 \quad \text{in } M.
\end{equation}
To obtain uniqueness for problem \eqref{e1}, we require that, for some \( 1 < p < 2 \),
\begin{equation}\label{e12}
\int_0^T\int_M |u(x,t)|^p \exp\{- h(x)\} \varphi^{2-p}(x) \, d\nu(x) dt < \infty.
\end{equation}
Clearly, if
\[
\varphi(x)=\exp\{- a[r(x)]^\gamma\} \quad \text{ for any } x\in M,
\]
for some $a>0$ and $\gamma>2$, then condition \eqref{e12} is substantially weaker than \eqref{ei6}. Obviously, in our approach it is essential to consider $p\neq 2$.

We will exhibit a class of potentials \( V \) for which the preceeding condition is satisfied, assuming that $M$ is a manifold with a pole. For this purpose, we study equation \eqref{ei2} in detail. Although a vast literature exists on this topic (see, e.g., \cite{Agm1, Agm2, Agm3, Agm4, Grig5, HPR, Pinch, Pinch2}), the existing results are not entirely sufficient for our goals (see Remarks \ref{decaysol}, \ref{oss-soprasol}, \ref{oss-soprasol2}). For specific potentials \( V \), on manifolds with a pole, we construct a solution with controlled pointwise decay at infinity. For others, we may only construct a positive locally Lipschitz weak supersolution to \eqref{ei2}, which still suffices to obtain a weaker uniqueness result for problem \eqref{e1}.

\smallskip

Our method is based on a suitable change of measure and the use of the associated weighted Laplacian. In fact, solutions $u$ to problem \eqref{e1} transform into solutions $w:=\frac{u}{\varphi}$ of problem \eqref{e35} whenever the
function \( \varphi \), which is also used to define the new measure \( d\mu \), is a solution
to the stationary Schr\"odinger equation \eqref{ei2}.  Therefore, by means of the uniqueness results recalled above for problem \eqref{e35}, we deduce that $w\equiv 0$, and thus $u\equiv 0$. In this regard, we note that condition \eqref{e12} is equivalent to \eqref{ei7}.

In some cases, however, we must rely on a function \( \varphi \) that is only a
weak supersolution of \eqref{ei2}. Despite this, we prove that any subsolution
to problem \eqref{e1} is mapped into a subsolution of problem \eqref{e35}.  As a consequence, we can establish the uniqueness of \emph{nonnegative} subsolutions
to \eqref{e1} that satisfy condition \eqref{e12}.

\bigskip

The paper is organized as follows. In Section~\ref{mb}, we recall preliminary notions from Riemannian geometry. Section~\ref{mr} states our main results, for both general and special classes of potentials \( V \). In Section~\ref{tow}, we explain how to reduce problem \eqref{e1} to \eqref{e35}. General uniqueness criteria are proved in Section \ref{proofgeneral}, while Section \ref{proofcs} addresses the cases of special potentials. Finally, in Section \ref{ex}, we describe some specific examples in which our uniqueness results are applied.

\section{Mathematical background}\label{mb}\setcounter{equation}{0}
Let $M$ be a connected Riemannian manifold equipped with a Riemannian metric $g$. For any smooth function $u$ defined on $M$, the gradient $\nabla u$ is a vector field on $M$. In local coordinates $x^1, \dots, x^n$, it takes the form
\[
(\nabla u)^i = g^{ij} \frac{\partial u}{\partial x^j},
\]
where we adopt the Einstein summation convention over repeated indices.

Now consider a smooth vector field $F$ on $M$. The divergence of $F$, denoted $\operatorname{div} F$, is a scalar function given in local coordinates by
\[
\operatorname{div} F = \frac{1}{\sqrt{\det g}} \frac{\partial}{\partial x^i} \left( \sqrt{\det g} F^i \right).
\]

The Riemannian volume form on $M$, denoted $\nu$, is expressed as
\[
d\nu = \sqrt{\det g} \, dx^1 \cdots dx^n.
\]

According to the divergence theorem, for any smooth function $u$ and any smooth vector field $F$ on $M$, with at least one of them having compact support, the following identity holds:
\[
\int_M u \, \operatorname{div} F \, d\nu = - \int_M \langle \nabla u, F \rangle \, d\nu,
\]
where $\langle \cdot, \cdot \rangle$ denotes the inner product induced by the metric $g$.

In particular, if $F = \nabla v$ for a smooth function $v$, we obtain
\begin{equation}\label{e2}
\int_M u \, \Delta v \, d\nu = - \int_M \langle \nabla u, \nabla v \rangle \, d\nu, \tag{2.2}
\end{equation}
provided that at least one of the functions $u$ or $v$ has compact support. Here, the operator
\[
\Delta := \operatorname{div} \circ \nabla
\]
is known as the Laplace (or Laplace-Beltrami) operator associated with the Riemannian manifold $M$.

From the identity above, we deduce the Green identities, which express the symmetry of the Laplace operator with respect to the $L^2$ inner product:
\begin{equation}\label{e3}
\int_M u \, \Delta v \, d\nu = - \int_M \langle \nabla u, \nabla v \rangle \, d\nu = \int_M v \, \Delta u \, d\nu.
\end{equation}

Let us consider a weighted manifold $(M, g, \mu)$, where the measure $\mu$ is defined via a positive, locally Lipschitz function $\varphi : M \to \mathbb{R}$ by
\[
d\mu = \varphi^2 \, d\nu,
\]
with $\nu$ the Riemannian volume measure on $M$.

The associated weighted divergence operator $\operatorname{div}_\mu$ is given by
\[
\operatorname{div}_\mu F := \frac{1}{\varphi^2} \operatorname{div}(\varphi^2 F),
\]
for any smooth vector field $F$ on $M$.

This leads to the definition of the weighted Laplace operator:
\begin{equation}
\Delta_\mu f := \operatorname{div}_\mu (\nabla f) = \frac{1}{\varphi^2} \operatorname{div}(\varphi^2 \nabla f),
\end{equation}
which, under the assumption that $\varphi$ is locally Lipschitz and strictly positive, can be expressed almost everywhere as
\[
\Delta_\mu f = \Delta f + \frac{2}{\varphi} \left\langle \nabla \varphi, \nabla f \right\rangle.
\]

Furthermore, the Green identity remains valid with respect to the weighted measure $\mu$, namely:
\begin{equation}
\int_M u \, \Delta_\mu v \, d\mu = - \int_M \langle \nabla u, \nabla v \rangle \, d\mu
= \int_M v \, \Delta_\mu u \, d\mu,
\end{equation}
provided that at least one of the functions \( u, v \) has compact support.

\subsection{Manifolds with a pole}
To study the behavior of the Laplacian of the distance function on a Riemannian manifold, let us fix a reference point \( x_0 \in M \), and denote by \( \mathrm{Cut}(x_0) \) the \emph{cut locus} of \( x_0 \). For any point \( x \in M \setminus \left( \mathrm{Cut}(x_0) \cup \{x_0\} \right) \), one can define \emph{polar coordinates} centered at \( x_0 \).

More precisely, for each such \( x \), there exists a unique minimizing geodesic joining \( x_0 \) to \( x \), whose initial direction defines an angle \( \theta \in \mathbb{S}^{n-1} \) in the unit tangent sphere \( T_{x_0}^1 M \). Identifying \( T_{x_0} M \cong \mathbb{R}^n \), we can regard \( \theta \) as a point on the standard Euclidean sphere \( \mathbb{S}^{n-1} \). The function
\[
r(x) := \operatorname{dist}(x, x_0)
\]
defines the \emph{radial coordinate} of \( x \), and the pair \( (r, \theta) \) yields a local coordinate system on \( M \setminus \left( \mathrm{Cut}(x_0) \cup \{x_0\} \right) \), referred to as \emph{polar coordinates} centered at \( x_0 \).

In these coordinates, the Riemannian metric can be expressed as
\[
ds^2 = dr^2 + A_{ij}(r, \theta) \, d\theta^i d\theta^j,
\]
where \( (\theta^1, \dots, \theta^{n-1}) \) are local coordinates on \( \mathbb{S}^{n-1} \), and the matrix \( (A_{ij}) \) is symmetric and positive definite. Denoting \( A := \det(A_{ij}) \), a straightforward computation shows that the Laplace--Beltrami operator takes the form
\begin{equation}\label{e4}
\Delta = \frac{\partial^2}{\partial r^2} + \mathfrak m(r, \theta) \frac{\partial}{\partial r} + \Delta_{S_r},
\end{equation}
where
\[
\mathfrak m(r, \theta) := \frac{\partial}{\partial r} \left( \log \sqrt{A(r, \theta)} \right),
\]
and \( \Delta_{S_r} \) denotes the Laplace--Beltrami operator on the geodesic sphere
\[
S_r := \partial B(x_0, r) \setminus \mathrm{Cut}(x_0).
\]

We say that \( M \) is a \emph{manifold with a pole} if there exists a point \( x_0 \in M \) such that \( \mathrm{Cut}(x_0) = \emptyset \). In this case, polar coordinates are globally defined on \( M \setminus \{x_0\} \). Moreover,
\[\mathfrak m(r, \theta)=\Delta r(x) \quad \text{ for any } x\in M\setminus\{x_0\}\,.\]

An important example of manifolds with a pole is provided by \emph{Cartan--Hadamard manifolds}, that is, complete, simply connected Riemannian manifolds with nonpositive sectional curvature. In fact, a fundamental property of such spaces is that the cut locus of any point \( x_0 \in M \) is empty. Moreover, on Cartan-Hadamard manifolds, Laplacian comparison theorems yield the lower bound
\begin{equation*}\label{e7}
\mathfrak m(r, \theta) \geq \frac{n - 1}{r} \quad \text{for all } x \equiv (r, \theta) \in M \setminus \{x_0\}.
\end{equation*}

\medskip

Let
\[
\mathcal{A} := \left\{ f \in C^\infty((0, \infty)) \cap C^1([0, \infty)) :
f(0) = 0,\ f'(0) = 1,\ f > 0 \text{ on } (0, \infty) \right\}.
\]
We say that a Riemannian manifold \( M \) is {\it spherically symmetric}, or a {\it model manifold}, if
its metric takes the form
\[
ds^2 = dr^2 + \psi(r)^2 d\theta^2,
\]
where \( d\theta^2 \) is the standard metric on the unit sphere \( \mathbb{S}^{n-1} \), and
\( \psi \in \mathcal{A} \). In such a case, we denote the manifold by \( M \equiv M_\psi \).

Moreover, the volume element has the form
\[
d\mu = \psi(r)^{n-1}dr d\theta.
\]
In addition, the Laplace--Beltrami operator reads
\begin{equation}\label{e4b}
\Delta = \frac{\partial^2}{\partial r^2}
+ (n - 1) \frac{\psi'(r)}{\psi(r)} \frac{\partial}{\partial r}
+ \frac{1}{\psi(r)^2} \Delta_{\mathbb{S}^{n-1}}.
\end{equation}

As notable examples:
 \( \psi(r) = r \) yields the Euclidean space \( \mathbb{R}^n \);
\( \psi(r) = \sinh r \) corresponds to the hyperbolic space \( \mathbb{H}^n \).

\subsection{Definition of solutions}
Whenever we refer to a solution, subsolution, or supersolution of the Cauchy problem or of the elliptic equation, we shall mean classical solutions, subsolutions, or supersolutions. We will also consider weak solutions, understood in the following sense.

\begin{definition}\label{defsol1}
A function \( w: M \times (0,T) \to \mathbb{R} \) is called a \emph{weak subsolution} of equation
\begin{equation}\label{ei1}
w_t - \Delta_\mu w = 0 \quad \text{in } M \times (0,T)\,.
\end{equation}
 if
\begin{itemize}
\item \( w \in L^2_{\mathrm{loc}}((0,T); H^1_{\mathrm{loc}}(M, d\mu)) \),
\item \( \partial_t w \in L^1_{\mathrm{loc}}((0,T); H^{-1}_{\mathrm{loc}}(M, d\mu)) \),
\item and for every nonnegative test function \( \psi \in C_c^\infty(M \times (0,T)) \), the following inequality holds:
\[
\iint_{M \times (0,T)} \left( -w \, \partial_t \psi + \langle \nabla w, \nabla \psi \rangle \right)\, d\mu dt \leq 0.
\]
\end{itemize}
\end{definition}

\begin{definition}\label{defsol2}
Let \( w_0 \in L^2_{\mathrm{loc}}(M, d\mu) \). A function
\[
w \in L^2_{\mathrm{loc}}((0,T); H^1_{\mathrm{loc}}(M, d\mu)) \cap C^0([0,T); L^2_{\mathrm{loc}}(M, d\mu))
\]
is a \emph{weak subsolution} of the Cauchy problem
\[
\begin{cases}
w_t - \Delta_\mu w \leq 0 & \text{in } M \times (0,T), \\
w(\cdot, 0) = w_0 & \text{in } M,
\end{cases}
\]
if
\begin{itemize}
\item \( w \) is a weak subsolution of equation \eqref{ei1},
\item and \( w(t) \to w_0 \) in \( L^2_{\mathrm{loc}}(M, d\mu) \) as \( t \to 0^+ \).
\end{itemize}
\end{definition}

\begin{definition}
Let \( \varphi \in H^1_{\mathrm{loc}}(M) \cap C^0(M) \). We say that \( \varphi \) is a \emph{weak supersolution} of equation \eqref{ei2}
if for all nonnegative test functions \( \eta \in C_c^\infty(M) \), the following inequality holds:
\begin{equation}\label{e30}
\int_M \langle \nabla \varphi, \nabla \eta \rangle\, d\nu + \int_M V(x) \varphi \eta \, d\nu \geq 0.
\end{equation}
We say that $\varphi$ is a \emph{weak subsolution} if $"\geq"$ is replaced by $"\leq"$.

\end{definition}

\section{Main results}\label{mr}\setcounter{equation}{0}

\subsection{General criteria}
Let $h:(0, \infty)\to \mathbb R$ be a function such that
\begin{equation}\label{e10}
 h\in C((0, \infty)), h>0, h \text{ is increasing}, \int_1^\infty \frac{r}{ h(r)}dr=\infty.
\end{equation}

Let $\mathcal H$ be $L^2(M)-$realization of the Schr\"odinger differential operator $-\Delta + V$, and let $\sigma (\mathcal H)$ denote its spectrum.
Set
\begin{equation}\label{e56}
\Lambda_1:=\inf \sigma(\mathcal H).
\end{equation}
  Note that $\Lambda_1\geq 0$ if and only if there exists a positive solution of equation \eqref{ei2} (see, e.g., \cite[Theorem 3.1, Remark 3,2]{Agm2} and \cite[Theorem 2.6]{Agm4}).

\begin{theorem}\label{t1} Let $M$ be a complete noncompact Riemannian manifold.
Let $u$ be a solution of problem \eqref{e1}. Let $ h$ be a function fulfilling hypothesis \eqref{e10}.
Assume that $\Lambda_1\geq 0$, and let $\varphi$ be a positive solution of equation \eqref{ei2}.
If, for some $1<p\leq2,$ condition \eqref{e12} is satisfied, then
\[u\equiv 0 \quad \text{ in } M\times (0, T)\,.\]
\end{theorem}

From Theorem it immediately follows the next
\begin{corollary}\label{cor1}  Let $M$ be a complete noncompact Riemannian manifold. Let $ h$ be a function fulfilling hypothesis \eqref{e10}.
Assume that $\Lambda_1\geq 0$, and let $\varphi$ be a positive solution of equation \eqref{ei2}.
Let $f\in C(M\times (0,T))$ and $u_0\in C(M)$. Then there exists at most one solution of problem 
\begin{equation}\label{e70}
\begin{cases}
u_t = \Delta u + f & \text{ in } M\times (0, T]\\
u= u_0 & \text{ in } M\times \{0\}\,.
\end{cases}
\end{equation}
which fulfills condition \eqref{e12}  for some $1<p<2$.
\end{corollary}

A similar corollariy can also be stated for Theorem \ref{t3b} below; we leave the formulation to the reader. 
Observe that in Example \ref{exRn} we discuss, in a particular case, not only uniqueness but also existence for problem \eqref{e70} with $u_0$ growing at infinity like $\exp\{K|x|^{\frac{\alpha}2+1}\}$
for some $\alpha>0, K>0$.

\begin{theorem}\label{t1a} Let $M$ be a complete noncompact Riemannian manifold.
Let $u$ be a nonnegative subsolution of problem \eqref{e1}. Let $ h$ be a function fulfilling hypothesis \eqref{e10}.
Assume that there exists a positive locally Lipschitz weak supersolution $\xi$ of equation \eqref{ei2}. If, for some $1<p\leq2,$
\begin{equation}\label{e12a}
\int_0^T\int_M |u(x,t)|^p \exp\{- h(x)\}\xi^{2-p}(x) d\nu(x) dt <\infty.
\end{equation}
then
\[u\equiv 0 \quad \text{ in } M\times (0, T)\,.\]
\end{theorem}

\begin{remark}{\em 
We assume that \( \xi \) is locally Lipschitz in order to define the measure \( d\mu := \xi^2\, d\nu \), which will be used to construct the weighted Laplacian \( \Delta_\mu \).
Under this regularity assumption, \( \Delta_\mu \) is well-defined almost everywhere in \( M \) (see Section \ref{mb}).} 
\end{remark}

\begin{remark}\label{Vuniq}
{\em
Note that in Theorem \ref{t1}  if
\begin{equation}\label{e40}
\varphi(x)=\exp\{- a[r(x)]^\gamma\} \quad \text{ for any } x\in M,
\end{equation}
for some $a>0$ and $\gamma>2$, then condition \eqref{e12} is much weaker than \eqref{ei6}, that is the condition usually requested in the case with $V\equiv 0$.

The same holds for Theorem \ref{t1a} with condition \eqref{e12} replaced by \eqref{e12a}, whenever
\[\xi(x)=\exp\{- a[r(x)]^\gamma\} \quad \text{ for any } x\in M,\]
for some $a>0$ and $\gamma>2$.}
\end{remark}

\subsection{Special cases}
In this Section we consider special classes of Riemannian manifolds. Specifically, we assume that
\begin{equation}\label{e5}
\begin{cases}
(i) &  M \mbox{ is a manifold with a pole } x_0\in M, \\
(ii) &  \mathfrak m(r, \theta)\geq 0  \mbox{ for any } x=(r, \theta)\in M\setminus\{x_0\}\,.
\end{cases}
\end{equation}

\smallskip

Now, we deal with potentials $V$ that satisfy the next assumption
\begin{equation}\label{e16b}
\begin{cases}
(i) &  V\in C^1(M) \\
(ii)  & V(x) \geq c_0 [r(x)]^{\alpha}  \mbox{ for all } x\in M\setminus B_{\bar R}, \text{ for some } \bar R>0, c_0>0, \alpha>0,
\end{cases}
\end{equation}
and we assume that
\[\Lambda_1 = 0\,.\]

Let
\begin{equation}\label{e55}
a_0:=\frac{2\sqrt{c_0}}{\alpha+2}\,.
\end{equation}

We prove the following
\begin{theorem}\label{t3b}
Let assumptions \eqref{e10}, \eqref{e5}, \eqref{e16b}  be satisfied. Suppose that $\Lambda_1=0$. Let $u$ be a solution of problem \eqref{e1}.
If, for some $1<p<2,$
\begin{equation}\label{e14}
\int_0^T\int_M |u(x,t)|^p \exp\{- h(x)\}\exp\left\{- a_0 (2-p) [r(x)]^{\left(\frac{\alpha}2+1\right)}\right\} d\nu(x) dt <\infty,
\end{equation}
then
\[u\equiv 0 \quad \text{ in } M\times (0, T)\,.\]
\end{theorem}

\begin{remark}\label{exe} Let us comment about the validity of \eqref{e16b} combined with $\Lambda_1=0.$
Consider a function $F\in C^1(M)$ fulfilling hypothesis \eqref{e16b}. Then, since $\displaystyle \lim_{r(x)\to +\infty} F(x)=+\infty$, the operator $-\Delta +F $ has compact resolvent, then its spectrum is purely discrete and the eigenvalues form a sequence which diverges to $+\infty$ (see, e.g., \cite[Theorem XIII.16]{RS}, where the equation posed in $\mathbb R^n$ is dealt with, but exactly the same holds on $M$). Let $\lambda_0$ be the first eigenvalue of $-\Delta +F $. Define
\[V:= F -\lambda_0.\]
Then assumption \eqref{e16b} is fulfilled and $\Lambda_1=0.$
\end{remark}

\bigskip

Let $\lambda_1$ be the infimum of the spectrum of $-\Delta$ in $L^2(M)$. Clearly, $\lambda_1\geq 0$.  We address potentials $V$ that fulfill the following hypothesis
\begin{equation}\label{e16}
\begin{cases}
(i) &  V\in C^1(M) \\
(ii) &  V(x) \geq -\lambda_1  \mbox{ for any  } x\in M,  \\
(iii)  & V(x) \geq c_0 [r(x)]^{\alpha}  \mbox{ for all } x\in M\setminus B_{\bar R}, \text{ for some } \bar R>0, c_0>0, \alpha>0.
\end{cases}
\end{equation}
In this case, is it easily seen that  (see Lemma \ref{lemma1})
\[\Lambda_1\geq 0\,.\]

For such potentials, we can prove the following
\begin{theorem}\label{t3}
Let assumptions \eqref{e10}, \eqref{e5}, \eqref{e16}  be satisfied. Let $u$ be a {\em nonnegative} subsolution of problem \eqref{e1}.
If, for some $1<p<2,$ \eqref{e14} holds, then
\[u\equiv 0 \quad \text{ in } M\times (0, T)\,.\]
\end{theorem}

\bigskip

\begin{remark}
In Theorems \ref{t3b} and \ref{t3}, a typical choice is $h(r)=r^2$. Note that if $\alpha>2$, then
\[\frac{\alpha}2+1>2.\]
Therefore \eqref{e14} is much weaker than \eqref{ei6}, that is the condition usually requested in the case with $V\equiv 0$.
\end{remark}

\section{Transformation to the Weighted Heat Equation}\label{tow}\setcounter{equation}{0}
\begin{proposition}\label{prop1}
Let $(M, g)$ be a complete, noncompact Riemannian manifold with Riemannian volume measure $d\nu$. Let $V \in C^1(M).$ Let $u$ be a classical solution to
\begin{equation}\label{eq:heat-potential}
\partial_t u - \Delta u + V(x) u = 0 \quad \text{in } M\times (0, T).
\end{equation}
Assume that $\varphi \in C^2(M)$ is a positive solution to \eqref{ei2}.
Define the function
\begin{equation}\label{e33}
w := \frac{u}{\varphi} \quad \text{in } M\times (0, T).
\end{equation}
Let $d\mu := \varphi^2 \, d\nu$.
Then $w$ satisfies the weighted heat equation \eqref{ei1}.
\end{proposition}

\noindent{\it Proof of Proposition \ref{prop1}.} Let $u$ be a classical solution to \eqref{eq:heat-potential}, and define $w := \frac{u}{\varphi}$. Using the product rule, we compute the gradient of $u$:
\[
\nabla u = \nabla(\varphi w) = w \nabla \varphi + \varphi \nabla w.
\]
Applying the Laplace--Beltrami operator, we obtain
\[
\Delta u = \Delta(\varphi w) = w \Delta \varphi + 2 \langle \nabla \varphi, \nabla w \rangle + \varphi \Delta w.
\]
Substituting this into equation \eqref{eq:heat-potential} gives
\[
\partial_t (\varphi w) - \left[ w \Delta \varphi + 2 \langle \nabla \varphi, \nabla w \rangle + \varphi \Delta w \right] + V \varphi w = 0.
\]
Expanding the time derivative and using the fact that $\partial_t \varphi = 0$, we find
\[
\varphi \partial_t w - w \Delta \varphi - 2 \langle \nabla \varphi, \nabla w \rangle - \varphi \Delta w + V \varphi w = 0.
\]
Divide by $\varphi$:
\[
\partial_t w - \Delta w - 2 \left\langle \frac{\nabla \varphi}{\varphi}, \nabla w \right\rangle + \left( -\frac{\Delta \varphi}{\varphi} + V \right) w = 0.
\]
Now, since $\varphi$ satisfies the stationary equation \eqref{ei2}, we have \(\Delta \varphi = V \varphi\), hence
\[
-\frac{\Delta \varphi}{\varphi} + V = 0.
\]
Thus the previous equation simplifies to
\begin{equation}\label{e32}
\partial_t w - \Delta w - 2 \left\langle \frac{\nabla \varphi}{\varphi}, \nabla w \right\rangle = 0.
\end{equation}

Recall that the weighted Laplacian with respect to the measure $d\mu = \varphi^2 d\nu$ is
\[
\Delta_\mu w = \Delta w + 2 \left\langle \frac{\nabla \varphi}{\varphi}, \nabla w \right\rangle.
\]
Therefore, \eqref{e32} yields
\[
\partial_t w - \Delta_\mu w = 0.\,  \hfill \square
\]

\begin{proposition}\label{prop2}
Let $(M, g)$ be a complete, noncompact Riemannian manifold with Riemannian volume measure $d\nu$. Let $V \in C^1(M).$  Suppose that
\begin{itemize}
  \item \( u \) is a (classical) nonnegative subsolution of \eqref{e1};
    \item \( \varphi \) is a positive weak supersolution of the stationary equation \eqref{ei2}.
  \end{itemize}
Define $w$ as in \eqref{e33}, and consider the weighted measure \( d\mu := \varphi^2 \, d\nu \). Then $w$ is a weak subsolution of the weighted heat equation
\eqref{ei1}.
\end{proposition}

\noindent{\it Proof of Proposition \ref{prop2}.} We aim to prove that for all nonnegative test functions \( \psi \in C_c^\infty(M \times (0,T)) \),
\[
\int_0^T\int_M \left( -w \partial_t \psi + \langle \nabla w, \nabla \psi \rangle \right)\, d\mu dt \leq 0.
\]
Let \( \psi \in C_c^\infty(M \times (0,T)) \), \( \psi \geq 0 \), and define \( \phi := \varphi \psi \). Then
\[
\phi_t = \varphi \psi_t, \quad \nabla \phi = \psi \nabla \varphi + \varphi \nabla \psi.
\]
Since \( u = \varphi w \), we have \( \nabla u = \varphi \nabla w + w \nabla \varphi \). The classical differential inequality for \( u \) gives
\[
\int_0^T\int_M \left( -u \phi_t + \langle \nabla u, \nabla \phi \rangle + V u \phi \right)\, d\nu dt \leq 0.
\]
Each term expands as
\begin{align*}
-u \phi_t &= -\varphi^2 w \psi_t, \\
\langle \nabla u, \nabla \phi \rangle &= \varphi^2 \langle \nabla w, \nabla \psi \rangle + \varphi \psi \langle \nabla w, \nabla \varphi \rangle + w \psi |\nabla \varphi|^2 + w \varphi \langle \nabla \varphi, \nabla \psi \rangle, \\
V u \phi &= V \varphi^2 w \psi.
\end{align*}

Substituting into the inequality, we obtain
\[
\begin{aligned}
\iint \Big[
& -\varphi^2 w \psi_t
+ \varphi^2 \langle \nabla w, \nabla \psi \rangle
+ \varphi \psi \langle \nabla w, \nabla \varphi \rangle \\
& + w \psi |\nabla \varphi|^2
+ w \varphi \langle \nabla \varphi, \nabla \psi \rangle
+ V \varphi^2 w \psi
\Big]\, d\nu dt \leq 0.
\end{aligned}
\]
Hence
\begin{equation}\label{e34}
\int_0^T\int_M \left( -w \psi_t + \langle \nabla w, \nabla \psi \rangle \right)\, d\mu dt\leq - \mathfrak I,
\end{equation}
where
\[\mathfrak I:=
\iint \left(
\psi \langle \nabla w, \nabla \varphi \rangle
+ w \psi \frac{|\nabla \varphi|^2}{\varphi^2}
+ w \langle \nabla \varphi, \nabla \psi \rangle
+ V w \psi
\right)\, d\mu dt.
\]
Recall that for every \( \eta \in C_c^\infty(M) \), \( \eta \geq 0 \), we have
\[
\int_M \langle \nabla \varphi, \nabla \eta \rangle + \int_M V \varphi \eta \geq 0.
\]
For each $t>0$, setting \( \eta := w \psi \), we comupte
\[
\langle \nabla \varphi, \nabla (w \psi) \rangle = \psi \langle \nabla \varphi, \nabla w \rangle + w \langle \nabla \varphi, \nabla \psi \rangle.
\]
Thus, multiplying by \( \varphi \) and  integrating in time over $(0,T)$ this matches:

\begin{equation}\label{e37}
\int_0^T \int_M \left(
\psi \langle \nabla w, \nabla \varphi \rangle
+ w \langle \nabla \varphi, \nabla \psi \rangle
+ V w \psi
\right)\, d\mu dt \geq 0.
\end{equation}

Since the extra term
\[
\int_0^T\int_M w \psi  \frac{|\nabla \varphi|^2}{\varphi^2}\, d\mu dt\geq 0,
\]
\eqref{e37} implies that $\mathfrak I\geq 0$. Thus, from \eqref{e34} the thesis follows.

$\hfill \square $

\section{Proofs of Theorems \ref{t1} and \ref{t1a}}\label{proofgeneral}\setcounter{equation}{0}
We need to recall the next uniqueness result for problem \eqref{e35} (see \cite{MeRo}; see also \cite[Theorem 11.9]{Grig3}, \cite[Theorem 2.2]{Pu15}).
\begin{theorem}\label{tlett1}
Let $M$ a complete noncompact RIemannian manifold.
Let $w$ be a solution of problem \eqref{e35}. Let $ h$ be a function fulfilling hypothesis \eqref{e10}.
If, for some $1<p\leq2$,
\begin{equation}\label{e36}
\int_0^T\int_M |u(x,t)|^p \exp\{- h(x)\} d\mu(x) dt <\infty,
\end{equation}
then
\[u\equiv 0 \quad \text{ in } M\times (0, T)\,.\]
\end{theorem}

\noindent{\it Proof of Theorem \ref{t1}.} Let $\varphi$ be a solution of equation \eqref{ei2}. Set $w:=\frac{u}{\varphi}$. In view of Proposition \ref{prop1}, $w$ satisfies \eqref{e35}, where $d\mu=\varphi^2 d\nu.$
In view of \eqref{e12}, since $u=\varphi w$, we have that
\[\int_0^T\int_M |u(x,t)|^p \exp\{-h(x)\}\varphi^{2-p}(x) d\nu(x) = \int_0^T\int_M |w(x,t)|^p \exp\{-h(x)\}d\mu(x) <\infty,\]
that is \eqref{e36}. By Theorem \ref{tlett1}, we can conclude that $w\equiv 0$. This implies that $u\equiv 0.$ $\hfill \square$

\bigskip

By minor changes in the proof of Theorem \ref{tlett1}, it possible to show the next result.
\begin{theorem}\label{tlett2}
Let $M$ a complete noncompact RIemannian manifold.
Let $w$ be a nonnegative weak subsolution of problem
\eqref{e35}. Let $h$ be a function fulfilling hypothesis \eqref{e10}.
If, for some $1<p\leq2,$ condition \eqref{e36} is satisfied,
then
\[u\equiv 0 \quad \text{ in } M\times (0, T)\,.\]
\end{theorem}

\noindent{\it Proof of Theorem \ref{t1a}} Let $\varphi$ be a weak supersolution of equation \eqref{ei2}. Set $w:=\frac{u}{\xi}$. In view of Proposition \ref{prop2}, $w$ is nonnegative weak subsolution to \eqref{e35}, where $d\mu=\xi^2 d\nu.$
In view of \eqref{e12}, since $u=\xi w$, we have that
\[\int_0^T\int_M |u(x,t)|^p \exp\{-h(x)\}\xi^{2-p}(x) d\nu(x) = \int_0^T\int_M |w(x,t)|^p \exp\{-h(x)\}d\mu(x) <\infty,\]
that is \eqref{e36}. By Theorem \ref{tlett2}, we can conclude that $w\equiv 0$. This implies that $u\equiv 0.$ $\hfill \square$

\section{Proofs of Theorems \ref{t3b} and \ref{t3}}\label{proofcs}\setcounter{equation}{0}

\subsection{Proof of Theorem \ref{t3b}}

Concerning the existence of positive solutions to equation \eqref{ei2} we recall the next result
(see, e.g., \cite[Theorem 3.1, Remark 3,2]{Agm2} and \cite[Theorem 2.6]{Agm4}).

\begin{proposition}\label{prop3}
Let $M$ be a complete noncompact Riemannian manifold, let $V\in C^1(M)$. We have that $\Lambda_1\geq 0$ if and only if there exists a positive solution $\eta$ of equation
\eqref{ei2}.
\end{proposition}

\begin{remark}\label{decaysol}{\em
In \cite[Section 6]{Agm2} some properties of solutions to \eqref{ei2} are stated. However, for a positive solution of \eqref{ei2}, we need a certain point-wise decay at infinity, that, on Riemannian manifolds, is not available in the literature. See also Remarks \ref{oss-soprasol} and \ref{oss-soprasol2} for further comments.}
\end{remark}

We show the following

\begin{proposition}\label{prop4}
Let hypotheses \eqref{e5} and \eqref{e16b} be satisfied. Then
\[z(x):=\exp\{- a_0 r^{\frac{\alpha}2+1}(x)\}, \quad x\in M\setminus B_{\bar R},\]
with $a_0$ defined in \eqref{e55}, is a supersolution to equation
\begin{equation}\label{e29}
\Delta z - V z = 0 \quad \text{ in } M\setminus B_{\bar R}\,.
\end{equation}
\end{proposition}

\noindent{\it Proof of Proposition \ref{prop4}.}
Let
\[z(x):=\exp\{- a r^\gamma(x)\}, \quad x\in M\setminus B_{\bar R},\]
where $a>0, \gamma>0$ are parameters to be chosen in the sequel. We set
\[z(x)\equiv z(r) \quad \text{ with } r\equiv r(x)\,.\]
We have
\[z'(r)=-a\gamma r^{\gamma-1}\exp{- a r^\gamma},\]
\[z''(r)=-a\gamma(\gamma-1)r^{\gamma-2}\exp{- a r^\gamma}+a^2\gamma^2r^{2\gamma-2}\exp{- a r^\gamma}\,.\]
Therefore, from \eqref{e4} and \eqref{e16b} we obtain that, for any $r>\bar R,$
\begin{equation}\label{e21}
\begin{aligned}
&\Delta z(r)- V(r) z(r)\\&\leq z(r)\left\{a^2\gamma^2r^{2\gamma-2}-a \gamma(\gamma-1)r^{\gamma-2}-\mathfrak m(r, \theta)a \gamma r^{\gamma-1}-c_0 r^\alpha\right\}\\
&\leq z(r)\left\{a^2\gamma^2r^{2\gamma-2}- c_0 r^\alpha\right\}\leq 0,
\end{aligned}
\end{equation}
provided that $a>0$ is small enough and
$\gamma\leq \frac{\alpha}{2}+1.$ From now on, we take $$\gamma=\frac{\alpha}{2}+1$$
and $$a=a_0,$$
with $a_0$ given by \eqref{e55}. For this choice, \eqref{e21} clearly holds. This completes the proof. $\hfill \square$

\begin{remark}\label{oss-soprasol}{\em
Let $\xi$ be a positive supersolution of equation \eqref{e29}. If we know a priori that a positive solution $\eta$ of \eqref{ei2} fulfills $\eta\in L^2(M),$ then we can apply \cite[Theorem 2.7 and Corollary 2.8]{Agm4}, which are stated in $\mathbb R^n$, but with the same
proof they also holds in $M$. Consequently, we can infer that there exists $\bar C>0$ such that
\[ \eta(x) \leq \bar C \xi(x) \quad \text{for any } x\in M\setminus B_{\bar R}\,.\]}
\end{remark}

Let $z$ be the supersolution to equation \eqref{e29} given by Proposition \ref{prop4}.

\begin{proposition}\label{prop10}
Let hypotheses \eqref{e5} and \eqref{e16b} be satisfied. Suppose that $\Lambda_1=0$. Then there exists a positive solution of equation \eqref{ei2} such that, for some $\bar C>0$,
\begin{equation}\label{e50}
\varphi(x)\leq \bar C z(x) \quad \text{ for all } x\in M\setminus B_{\bar R}\,.
\end{equation}
\end{proposition}
\noindent{\it Proof of Proposition \ref{prop10}.} Since $\displaystyle \lim_{r(x)\to +\infty} V(x)=+\infty$, the operator $-\Delta +V $ has compact resolvent, then its spectrum is purely discrete, and the eigenvalues form a sequence which diverges to $+\infty$  (see, e.g., \cite[Theorem XIII.16]{RS}). Due to the fact that $\Lambda_1=0$, the first eigenvalue is $0$, therefore, there exists a positive solution
$\varphi\in L^2(M)$ of equation \eqref{ei2}. Since $z$ is a positive supersolution of equation \eqref{e29}, by \cite[Corollary 2.8]{Agm4} (see also Remark \ref{oss-soprasol}), there exists $\bar C>0$ such that
\[\varphi(x)\leq \bar C z(x) \quad \text{ for any } x\in M\setminus B_{\bar R}\,.\]
Hence the thesis follows. $\hfill \square$

\medskip

\noindent{\it Proof of Theorem \ref{t3b}} The thesis easily follows by combining Theorem \ref{t1} and Proposition \ref{prop10}. $\hfill \square$

\subsection{Proof of Theorem \ref{t3}}
\begin{lemma}\label{lemma1}
Let $M$ be a complete noncompact Riemannian manifold. Let assumption \eqref{e16} be satisfied. Then
\begin{equation}\label{e66}
\Lambda_1 \geq 0.
\end{equation}
\end{lemma}
\noindent{\it Proof of Lemma \ref{lemma1}.} Consider the quadratic form
\[Q[\zeta]:=\int_M \left\{ |\nabla \zeta|^2+ V(x) \zeta^2 \right\} d\nu\,.\]
We have that
\[\Lambda_1=\inf \big\{Q[\zeta]\,:\, \zeta\in C^\infty_c(M), \|\zeta\|_{L^2(M)}=1 \big\}\,.\]
Let $\zeta\in C^\infty_c(M)$ with  $\|\zeta\|_{L^2(M)}=1$.  Obviously, in view of the definition of $\lambda_1$,
\begin{equation}\label{e67}
\int_M |\nabla \zeta|^2 d\nu \geq \lambda_1 \int_M \zeta^2 d\nu\,.
\end{equation}
Then, in view of \eqref{e67} and hypothesis \eqref{e16},
\[Q[\varphi]\geq \int_M \left\{ |\nabla \zeta|^2-  \lambda_1 \zeta^2 \right\} \geq 0\,. \]
Hence \eqref{e66} follows. $\hfill \square$

\begin{proposition}\label{prop5}
Let hypotheses \eqref{e5} and \eqref{e16} be satisfied.  Then there exists a weak positive supersolution $\xi$ of equation \eqref{ei2}
which is locally Lipschitz in $M$ and satisfies
\begin{equation}\label{e18}
\xi(x) \leq \exp\{-a_0 r(x)^{\frac{\alpha}{2}+1}\} \quad \text{ for all } x\in M\setminus B_{\bar R}(x_0),
\end{equation}
with $a_0$ given by \eqref{e55}.
\end{proposition}

\noindent{\it Proof of Proposition \ref{prop5}.}
Consider the differential operator $\mathcal K:= -\Delta -\lambda_1$. Clearly, the infimum of the $L^2(M)-$spectrum of $\mathcal K$ is $0$. Therefore, we can apply Proposition \ref{prop3} to infer that there exists
 a positive solution $\eta$ of \eqref{ei2} with $V\equiv -\lambda_1$. Thus
\[\Delta \eta + \lambda_1 \eta = 0 \quad \text{ in } M\,. \]
In view of assumption \eqref{e16},
\[\Delta \eta - V(x) \eta \leq \Delta \eta + \lambda_1 \eta =0  \quad \text{ in } M\,.\]
Therefore $\eta$ is a supersolution of equation \eqref{ei2}.  Let $z$ be the supersolution given by Proposition \ref{prop4}.
Consider any $R_1, R_2$ so that $0<R_1<\bar R< R_2$. We can select $\underline C>0$ such that
\begin{equation}\label{e22}
\underline C\eta \leq z \quad \text{ in } B_{R_2}\setminus B_{R_1}\,
\end{equation}
Set
\[\tilde \zeta:=\min\{\underline C \eta, z\}\,.\]
It is easily seen that $\tilde \zeta$ is a weak supersolution of equation \eqref{e29}; moreover, $\tilde \zeta$ is locally Lipschitz in $M$. Define
\[ \xi:=\begin{cases}
    \tilde \zeta, & \mbox{in } M\setminus B_{\bar R} \\
    \underline C\eta, & \mbox{in } B_{\bar R}.
  \end{cases}\]
Thanks to \eqref{e22} it is direct to see that $\xi$ is a weak supersolution of equation \eqref{ei2}, and it is locally Lipschitz in $M$.
Furthermore, \eqref{e18} holds. $\hfill \square$


\begin{remark}\label{oss-soprasol2}{\em
(i) Observe that in view of \cite[Theorem 3.1]{Agm2}, there exists a positive solution $\varphi $ of equation \eqref{ei2} if and only if there exists a weak supersolution $\overline{ \varphi}$ of equation \eqref{ei2}. However, the proof is not constructive, therefore in general it is not guaranteed that $\varphi \leq C \overline{\varphi}$ in $M$, for some $C>0$. Therefore, we cannot use that result combined with Proposition \ref{prop5} to infer the existence of a global solution of
\eqref{ei2} with a certain point-wise decay at infinity.

\smallskip

\noindent (ii) Note that, under the assumptions of Theorem \ref{t3}, we have $\Lambda_1 > 0$.
Therefore, by Proposition \ref{prop3}, a positive solution $\eta$ of \eqref{ei2} exists.
Moreover, thanks to Proposition \ref{prop5}, there also exists a weak positive supersolution
$\xi$ of \eqref{ei2}.  Applying \cite[Theorem 4.2]{Agm2}, we deduce the existence of a constant $C > 0$,
a \textit{non-empty} compact set $K \subset M$, and a positive solution
$\overline{v}$ of the equation
\[
\Delta \overline{v} - V \overline{v} = 0 \quad \text{in } M \setminus K,
\]
such that
\[
\overline{v} \leq C \xi \quad \text{in } M \setminus K.
\]
We emphasize that in \cite[Theorem 4.2]{Agm2}, the compact set $K$ must be
non-empty; otherwise, simple counterexamples can be constructed.
As a consequence, by this method, we cannot obtain a \textit{global} positive solution to
equation \eqref{ei2} that decays at most like $\xi$ at infinity.}
\end{remark}

\begin{remark}\label{methodsoprasol}
{\em
In order to extend the argument used in the proof of Theorem \ref{t3b}
to the case of potentials satisfying condition \eqref{e16},
a global positive solution to equation \eqref{ei2} would be required. However, since we are unable to construct such a global solution
with a precise decay rate at infinity (see also Remark \ref{oss-soprasol2}),
we work directly with the weak supersolution $\xi$. As a consequence, the function \( w := \frac{u}{\xi} \)
is only a weak subsolution to problem \eqref{e35}.
}

\end{remark}

\noindent{\it Proof of Theorem \ref{t3}.} The thesis easily follows by combining Theorem \ref{t1a} and Proposition \ref{prop5}. $\hfill \square$

\section{Examples}\label{ex}

Now we discuss three examples in which, keeping the notation of Section \ref{mr}, \eqref{e40} is satisfied. Then $V$ influences the integral uniqueness condition for solutions of problem \eqref{e1} (see Remark \ref{Vuniq}).

\begin{example}\label{exRn} {\em Let $M=\mathbb R^n, n\geq 2.$
Consider the radial potential
\[
V(x) = -ab(n + b - 2) |x|^{b - 2} + a^2 b^2 |x|^{2b - 2},
\]
where \( a > 0 \), \( b > 2 \), and \( x \in \mathbb{R}^n \).  Let us analyze the sign of \( V(x) \). Write \( r = |x| \). Then
\[
V(r) = r^{b - 2} \left[ -ab(n + b - 2) + a^2 b^2 r^b \right].
\]
Define the critical radius
\[
r_0 := \left( \frac{n + b - 2}{ab} \right)^{1/b}.
\]
Then
\begin{itemize}
    \item \( V(r) < 0 \) for \( 0 < r < r_0 \),
    \item \( V(r) = 0 \) for \( r = r_0 \),
    \item \( V(r) > 0 \) for \( r > r_0 \).
\end{itemize}
Furthermore,  \( V(0) = 0 \). Thus, \( V(x) \) is a sign-changing potential: negative near the origin, vanishing at \( r = r_0 \), and positive for large \( |x| \).

\smallskip

\noindent (a) {\it Uniqueness for problem \eqref{e1}}.

We aim to construct a positive radial function \( \varphi(x) \) solving
\[
\Delta \varphi(x) - V(x) \varphi(x) = 0 \quad \text{in } \mathbb{R}^n.
\]

We claim that the function
\[
\varphi(x) = e^{-a |x|^b}
\]
satisfies the above equation. Let \( r = |x| \). Since \( \varphi \) is radially symmetric, we compute the Laplacian using the radial formula:
\[
\Delta \varphi(x) = \varphi''(r) + \frac{n-1}{r} \varphi'(r),
\]
where
\[
\varphi'(r) = -ab r^{b - 1} e^{-a r^b}, \quad
\varphi''(r) = \left[ -ab(b - 1) r^{b - 2} + a^2 b^2 r^{2b - 2} \right] e^{-a r^b}.
\]
Substituting into the expression for \( \Delta u \), we obtain
\[
\Delta \varphi(x) = e^{-a r^b} \left[ -ab(n + b - 2) r^{b - 2} + a^2 b^2 r^{2b - 2} \right],
\]
which gives
\[
\Delta \varphi(x) = V(x) \varphi(x),
\]
as claimed. This provides an explicit example of a smooth, radially symmetric solution to a Schr\"odinger-type equation with a sign-changing potential. In addition, by Proposition \ref{prop3}, $\Lambda_1\geq 0.$

By Theorem \ref{t1}, problem \eqref{e1} admits at most a unique solution $u$, provided that, for some $1<p<2,$
\begin{equation}\label{e74}
\int_0^T\int_{\mathbb R^n} |u(x,t)|^p \exp\{- h(x)\} \exp\{- a(2-p) |x|^{b}\}dx dt <\infty, 
\end{equation}
with $h$ fulfilling \eqref{e10}.

\smallskip

\noindent (b) {\it Existence and uniqueness for problem \eqref{e70}}. In view of the choice of $V$ made in (a), we can write that 
\( V \in C(\mathbb{R}^n) \), and there exist constants \( \alpha > 0 \), \( R > 1, c_0>0 \), and \( M > 0 \) such that
\begin{align*}
V(x) &\geq c_0 |x|^\alpha \quad \text{for all } |x| > R, \\
V(x) &\geq -M \quad \text{for all } |x| \leq R.
\end{align*}
In particular, $\alpha=2b-2>2.$ Define
\[
\bar u(x,t) := \exp\left\{ A(1 - Qt) |x|^b \right\},
\]
where \( A > 0 \), \( Q > 0 \) are parameters to be appropriately chosen. 

We claim that there exist  \( A, Q > 0, T>0 \) such that \( \bar u \) is a supersolution to
\[
\bar u_t - \Delta \bar u + V(x) \bar u = 0 \quad \text{in } \mathbb{R}^n \times (0, T).
\]
Indeed, let \( r = |x| \) and define
\[
\theta(x,t) := A(1 + Qt) r^b, \quad \text{so that} \quad \bar u(x,t) = e^{\theta(x,t)}.
\]
We compute
\[
\theta_t = AQ r^b,
\quad
\nabla \theta = A(1 + Qt) \beta r^{b - 2} x,
\quad
|\nabla \theta|^2 = A^2(1 + Qt)^2 b^2 r^{2b - 2},
\]
\[
\Delta \theta = A(1 + Qt) b(n + b - 2) r^{b - 2}.
\]
Then
\[
\begin{aligned}
\bar u_t &= \theta_t e^\theta = AQ r^b e^\theta, \\
\Delta \bar u &= (\Delta \theta + |\nabla \theta|^2) e^\theta.
\end{aligned}
\]
So,
\[
\bar u_t - \Delta \bar u + V(x) \bar u =
\left( \theta_t - \Delta \theta - |\nabla \theta|^2 + V(x) \right) e^\theta.
\]
First, consider the region \( |x| > R \). Here \( V(x) \geq c_0 r^\alpha = c_0 r^{2b - 2} \). Then
\[
\begin{aligned}
\theta_t - \Delta \theta - |\nabla \theta|^2 + V(x)
&= AQ r^b
- A(1 + Qt) b(n + b - 2) r^{b - 2} \\
&\quad - A^2(1 + Qt)^2 b^2 r^{2b - 2}
+ c_0 r^{2b - 2}.
\end{aligned}
\]
Divide both sides by \( r^{2b - 2} \). Then in order to get 
\begin{equation}\label{e71}
\bar u_t - \Delta \bar u + V(x) \bar u \geq 0 \quad \text{ in } [\mathbb R^n\setminus B_R(0)]\times (0, T),
\end{equation}
we need that 
\[
AQ r^{-(b - 2)}
- A(1 + Qt) b(n + b - 2) r^{-b}
- A^2(1 + Qt)^2 b^2
+ c_0 \geq 0.
\]
So it suffices to ensure
\[
A^2 b^2 (1 + Qt)^2 \leq \frac{c_0}{2},
\quad
AQ R^{-(b - 2)} \leq \frac{c_0}{4},
\quad
A b(n + b - 2) R^{-b} \leq \frac{c_0}{4}.
\]
We may guarantee this on \( t \in [0, T] \) by choosing
\[
A \leq \frac{1}{b} \sqrt{\frac{c_0}{2}}, \quad
Q \leq \min\left\{
\frac{c_0}{4A} R^{-(b - 2)}, \,
\frac{c_0}{4A b(n + b - 2)} R^{-b}
\right\}, \quad
T < \frac{1}{2Q}.
\]
Now, consider the region \( |x| \leq R \). Here \( r \leq R \), and \( V(x) \geq -M \), for some \( M > 0 \).
We estimate
\[
\begin{aligned}
\theta_t - \Delta \theta - |\nabla \theta|^2 + V(x)
&\geq AQ R^b
- A(1 + Qt) b(n + b - 2) R^{b - 2} \\
&\quad - A^2(1 + Qt)^2 b^2 R^{2b - 2} - M.
\end{aligned}
\]
We fix \( A, Q \) as before, and define \( C := \sup_{t \in [0, T]} (1 + Qt)^2 \leq (1 + QT)^2 \).  
So, if  we require
\[
AQ R^b \geq A b(n + b - 2) R^{b- 2}
+ A^2 C b^2 R^{2b - 2} + M,
\]
then
\begin{equation}\label{e72}
\bar u_t - \Delta \bar u + V(x) \bar u \geq 0 \quad \text{ in } B_R(0)\times (0, T)\,.
\end{equation}
This inequality holds for small \( A \), since the left-hand side is linear in \( A Q \), and the right-hand side is \( O(A) + O(A^2) + M \). So choosing \( A \), \( Q \) small enough, 
\eqref{e72} is fulfilled. Putting together \eqref{e71} and \eqref{e72} we obtain the conclusion of the previous claim. 

\smallskip

Now, fix $1<p<2.$ We select possibly smaller $A>0, Q>0$ such that 
\begin{equation}\label{e73}
K:=A(1+QT)<a(2-p)\,.
\end{equation}
Consider problem \eqref{e70} with $f=0$ and 
\[u_0(x)=\bar u(x, 0) \quad \text{ for all } x\in \mathbb R^n\,.\]
 For any $j\in \mathbb N, j\geq 2,$
let $u_j$ be the solution of problem 
\begin{equation}\label{e75}
\begin{cases}
\partial_t u_j - \Delta u_j + V(x) u_j= 0 & \text{ in } B_j(0)\times (0, T)\\
u_j = 0 & \text{ in } \partial B_j(0)\times (0,T)\\
u_j = u_0 \chi_j & \text{ in } B_j(0)\times\{0\},
\end{cases}
\end{equation}
where $\chi_j \in C^\infty_c B_j(0), \chi_j \equiv 1$ in $B_{j-1}, 0\leq\chi_j\leq 1$ for each $j\in \mathbb N$. 
Therefore, $\underline u\equiv 0$ is a subsolution of problem \eqref{e50}, while $\bar u$  is a supersolution. 
Clearly, by the comparison principle, 
\[0\leq u_j \leq \bar u \quad \text{ in } B_j(0)\times (0, T)\,.\]
By standard a priori estimates and a diagonal procedure, there exists a subsequence $\{u_{j_k}\}\subset \{u_j\}$ which converges in $C^{2,1}_{\mathrm{loc}}(\mathbb R^n\times [0, T])$ to some funciton $u$, which is a solution to problem 
\eqref{e70}. Furthermore, 
\[0\leq u \leq \bar u \quad \text{ in } \mathbb R^n\times (0, T)\,.\]
Hence, in view of \eqref{e73},
\[
\begin{aligned}
&\int_0^T\int_{\mathbb R^n} |u(x,t)|^p \exp\{- h(x)\} \exp\{- a(2-p) |x|^{b}\}dx dt \\ & \leq \int_0^T \int_{\mathbb R^n} |\bar u(x,t)|^p \exp\{- h(x)\} \exp\{- a(2-p) |x|^{b}\}dx dt \\
&\leq T \int_{\mathbb R^n} \exp\{K p|x|^b \} \exp\{- h(x)\} \exp\{- a(2-p) |x|^{b}\}dx<\infty\,.
\end{aligned}
\]
By Corollary \ref{cor1}, $u$ is the unique solution to \eqref{e70} which fulfills condition \eqref{e74}.
}
\end{example}

\begin{example}{\em
Let \( M \) be a model manifold with metric
\[
ds^2 = dr^2 + \psi^2(r)\, d\theta^2,
\]
where \( r > 0 \), \( \theta \in \mathbb{S}^{n-1} \), and \( \psi : [0, \infty) \to [0, \infty) \) is a smooth warping function with \( \psi(0) = 0 \), \( \psi'(0) = 1 \), and \( \psi(r) > 0 \) for \( r > 0 \).

We consider the radial potential
\[
V(r) = -ab(b - 1) r^{b - 2} + a^2 b^2 r^{2b - 2} - ab(n - 1) \frac{\psi'(r)}{\psi(r)} r^{b - 1},
\]
where \( a > 0 \) and \( b > 2 \) are parameters. We now define the radial function
\[
\varphi(r) = e^{-a r^b}.
\]
We claim that \( \varphi \) satisfies the equation
\[
\Delta \varphi(r) - V(r) \varphi(r) = 0 \quad \text{on } M.
\]

Indeed, computing:
\begin{align*}
\varphi'(r) &= -ab r^{b - 1} e^{-a r^b}, \\
\varphi''(r) &= \left[ -ab(b - 1) r^{b - 2} + a^2 b^2 r^{2b - 2} \right] e^{-a r^b},
\end{align*}
and using \eqref{e4b}, we obtain
\[
\Delta \varphi(r) = \left[ -ab(b - 1) r^{b - 2} + a^2 b^2 r^{2b - 2} - ab(n - 1) \frac{\psi'(r)}{\psi(r)} r^{b - 1} \right] e^{-a r^b} = V(r) \varphi(r).
\]
By Theorem \ref{t1}, problem \eqref{e1} admits a unique solution $u$, provided that , for some $1<p<2,$
\[\int_0^T\int_{M} |u(x,t)|^p \exp\{- h(r)\} \exp\{- a(2-p) r^{b}\}\psi^{n-1}(r) dr d\theta dt <\infty, \]
with $h$ fulfilling \eqref{e10}. Here we have used the fact that $d\mu=\psi^{n-1}(r) dr d\theta$, where $d\theta$ is the volume element on the unit sphere $\mathbb S^{n-1}.$

\smallskip

Let us analyze the behavior of \( V(r) \) both near the origin and at infinity, further assuming that
\[ \psi(r) = e^{\alpha r^{\gamma}}  \quad \text{ for all } r>1\,,\]
for some $\alpha>0, \gamma\geq 1.$
We assume that \( \psi(r) = r + o(r) \) as \( r \to 0 \), which corresponds to the standard Euclidean behavior near the pole. In this case,
\[
\frac{\psi'(r)}{\psi(r)} \sim \frac{1}{r} \quad \text{as } r \to 0,
\]
so the potential behaves like:
\[
V(r) \sim -ab(n + b - 2) r^{b - 2}, \quad \text{as } r \to 0.
\]
In particular, \( V(r) \to 0 \). For large \( r \), \[
\frac{\psi'(r)}{\psi(r)} = \alpha \gamma r^{\gamma - 1}.
\]
Therefore, the potential becomes:
\[
V(r) = -ab(b - 1) r^{b - 2} + a^2 b^2 r^{2b - 2} - ab(n - 1)\alpha \gamma r^{b + \gamma - 2}.
\]
We compare the three terms as \( r \to \infty \). The dominant behavior depends on the relative sizes of the exponents:
the first term grows like \( r^{b - 2} \),  the second term grows like \( r^{2b - 2} \), the third term grows like \( r^{b + \gamma - 2} \). Since \( \gamma \geq 1 \) and \( b > 2 \), we always have:
\[
b + \gamma - 2 \leq 2b - 2,
\]
with equality only if \( \gamma = b \). Therefore:
\begin{itemize}
    \item If \( \gamma < b \), then the dominant term is \( a^2 b^2 r^{2b - 2} \), so
    \[
    V(r) \sim a^2 b^2 r^{2b - 2} \to +\infty \text{ as } r\to +\infty;
    \]
    \item If \( \gamma = b \), then the second and third terms are both of order \( r^{2b - 2} \), thus $V(r) \to \operatorname{sgn}(ab-(n-1)\alpha\gamma)\infty$  \text{ as } $r\to +\infty$;
    \item If \( \gamma > b \), then the third term dominates and is negative, so
    \[
    V(r) \sim -ab(n - 1)\alpha \gamma r^{b + \gamma - 2} \to -\infty  \text{ as } r\to +\infty.
    \]
\end{itemize}}

\end{example}

\begin{example}{\em
Let $M=\mathbb{H}^n$. Suppose that $V\in C^1(\mathbb H^n)$ with
\[V\geq - \frac{(n-1)^2}{4} \quad \text{ in } \mathbb H^n,\]
\[V \geq c_0 [r(x)]^{\frac{\alpha}2+1} \quad \text{for all } x\in \mathbb H^n, r(x) >1\]
for some $c_0>0, \alpha>0$.

Since $\lambda_1=\frac{(n-1)^2}{4} ,$ we are in position to apply Theorem \ref{t3} to infer that
problem \eqref{e1} admits at most a unique nonnegative subsolution $u$, provided that, for some $1<p<2,$
\[\int_0^T\int_{M} |u(x,t)|^p \exp\{- h(r)\} \exp\{- a_0 (2-p) r^{\frac{\alpha}2+1}\}\psi^{n-1}(r) dr d\theta dt <\infty, \]
with $h$ fulfilling \eqref{e10}.}
\end{example}

Now, we discuss two examples in which, keeping the notation of Section \ref{mr}, \eqref{e40} is not satisfied. Then $V$ does not substantially influence the integral uniqueness condition (see Remark \ref{Vuniq}) for solutions of problem \eqref{e1} .

\begin{example}{\em
Recall that the Green function \( G(x,y) \) of a manifold
\( (M, g) \) is defined by
\[
G(x,y) := \int_0^\infty p(x,y, t) \, dt,
\]
where \( p(x,y, t) \) is the heat kernel. 

A continuous function \( V \geq 0 \) in \( M \) is said to be
\emph{Green bounded} (see \cite[Section 10]{Grig2}) if either \(  V \equiv 0 \), or the Green function
\( G(x,y) \) is finite (that is, \( (M,g) \) is non-parabolic), and
\[
\sup_{x \in M} \int_M G(x,y)\, V(y) \, d\nu(y) < \infty.
\]
By \cite[Lemma 10.3]{Grig2}, if $V$ is Green bounded, then there exists a positive solution $\varphi$ of equation \eqref{ei2} such that
\begin{equation}\label{e1ex}
\delta\leq \varphi\leq 1 \quad \text{ in } M,
\end{equation}
for some $0<\delta<1$.  Note that, in particular, if $\operatorname{supp} V$ is compact, then $V$ is Green bounded.
By Theorem \ref{t1}, problem \eqref{e1} admits at most only one solution $u$, provided that, for some $1<p<2,$
\begin{equation}\label{e2ex}
\int_0^T\int_{\mathbb R^n} |u(x,t)|^p \exp\{- h(x)\} \varphi(x)^{(2-p)}\}d\mu(x) dt <\infty,
\end{equation}
with $h$ fulfilling \eqref{e10}. In view of \eqref{e1ex}, condition \eqref{e2ex} is equivalent to
\[\int_0^T\int_{\mathbb R^n} |u(x,t)|^p \exp\{- h(x)\} d\mu(x) dt <\infty,\]
which is the same as in the case with $V\equiv 0$.
}
\end{example}

\begin{example}{\em
Let $M=\mathbb R^n.$ Consider a potential $V\in C^1(\mathbb R^n)$
fulfilling
\[
V(x) = \frac{b}{|x|^2} \quad \text{for any } |x| > 1,
\]
where \( b>0\).

Let us look for a smooth, positive function \( \varphi \) on \( \mathbb{R}^n \) such that
\[
\Delta \varphi - V(x)\varphi = 0 \quad \text{in } \mathbb{R}^n.
\]
For \( |x| > 1 \), set
\[
\varphi(x) = |x|^\beta,
\]
with \( \beta \in \mathbb{R} \). A direct computation gives
\[
\frac{\Delta \varphi}{\varphi} = \frac{\beta^2 + (n - 2)\beta}{|x|^2},
\]
so that
\[
\Delta \varphi - V(x)\varphi = 0 \quad \text{for all } |x| > 1
\]
is satisfied if and only if
\[
\beta^2 + (n - 2)\beta = b.
\]
The solutions to this equation are given by
\[
\beta_\pm = \frac{-(n - 2) \pm \sqrt{(n - 2)^2 + 4b}}{2}.
\]
The smaller root \( \beta_- \) is negative for all \( b > 0 \), while the larger root \( \beta_+ \) is positive.
We set
\[\beta:=\beta_-\,.\]
Now, extend \( \varphi \) to a smooth, positive function on all of \( \mathbb{R}^n \) (e.g., using a cutoff or mollification that preserves positivity and regularity). Define
\[
V(x) := \frac{\Delta \varphi(x)}{\varphi(x)} \quad \text{for } |x| \leq 1.
\]
By construction, \( V \in C^\infty(\mathbb{R}^n), V>0 \), and the equation
\[
\Delta \varphi - V(x)\varphi = 0
\]
holds in all of \( \mathbb{R}^n \).
By Theorem \ref{t1}, problem \eqref{e1} admits at most only one solution $u$, provided that, for some $1<p<2,$
\begin{equation}\label{e3ex}
\int_0^T\int_{\mathbb R^n} |u(x,t)|^p \exp\{- h(x)\} |x|^{(2-p)\beta_-}\}dx dt <\infty,
\end{equation}
with $h$ fulfilling \eqref{e10}. Observe that \eqref{e40} is not verified. Obviously, if \eqref{e3ex} is satisfied, then  also
\begin{equation}\label{e4ex}
\int_0^T\int_{\mathbb R^n} |u(x,t)|^p \exp\{- \mathfrak h(x)\} dx dt <\infty,
\end{equation}
holds, for some $\mathfrak h$ as in \eqref{e10}.  Clearly, \eqref{e4ex} is the same condition as in the case with $V\equiv 0$.
}
\end{example}

\end{document}